\documentclass[10pt]{elsarticle}
\usepackage{amssymb}
\usepackage{amsfonts}
\usepackage{amsthm}
\usepackage{amsmath}
\usepackage[left=2.2cm,top=1.75cm,right=2.2cm]{geometry}
\newtheorem{lem}{Lemma}
\newtheorem{remark}{Remark}
\newtheorem{theorem}{Theorem}
\numberwithin{equation}{section}

\begin{document}
\journal{Elsevier }
\begin{frontmatter}
\title{On statistical approximation properties of $q$-Baskakov-Sz\'{a}sz-Stancu operators}

\author[label1,label**] { Vishnu Narayan Mishra}
\ead{vishnunarayanmishra@gmail.com, vishnu\_narayanmishra@yahoo.co.in, v\_n\_mishra\_hifi@yahoo.co.in}

\author[label1]{Preeti Sharma}
\ead{preeti.iitan@gmail.com}

\author[label2,label3] { Lakshmi Narayan Mishra}
\ead{lakshminarayanmishra04@gmail.com, l\_n\_mishra@yahoo.co.in}

\fntext[label**]{Corresponding author}

\address[label1]{Department of Applied Mathematics \& Humanities,
Sardar Vallabhbhai National Institute of Technology, Ichchhanath Mahadev Dumas Road, Surat - 395 007 (Gujarat), India}

\address[label2]{L. 1627 Awadh Puri Colony Beniganj, Phase -III, Opposite - Industrial Training Institute (I.T.I.), Ayodhya Main Road, Faizabad - 224 001, (Uttar Pradesh), India}

\address[label3]{Department of Mathematics, National Institute of Technology, Silchar - 788 010, Cachar (Assam), India}

\begin{abstract}
In the present paper, we consider  Stancu type generalization of Baskakov-Sz\'{a}sz operators based on the $q$-integers and obtain statistical and weighted statistical approximation properties of these operators. Rates of statistical convergence by means of the modulus of continuity and the Lipschitz type maximal function are also established for operators.
\end{abstract}

\end{frontmatter}

\section{Introduction}
In the recent years several operators of summation-integral type have been proposed and their approximation properties have been discussed. In the present paper our aim is to investigate statistical approximation properties of a Stancu type $q$-Baskakov-Sz\'{a}sz operators.
Firstly, Baskakov-Sz\'{a}sz operators based on $q$-integers was introduced by Gupta \cite{vj} and established some approximation results. The $q$-Baskakov-Sz\'{a}sz operators is defined as follows:

\begin{equation}\label{eq1}
\mathcal{D}_{n}^{q}(f,x)=  [n]_q \sum_{k=0}^{\infty}p^q_{n,k}(x)\int_{0}^{q/{1-q^n}}q^{-k-1}s^q_{n,k}(t)f\bigg(t \,q^{-k}\bigg) d_qt,
\end{equation}
where $x\in [0,\infty)$ and \begin{equation}
p^q_{n,k}(x)=\left[\begin{array}{c}n+k-1\\k \end{array}\right]_q q^{k(k-1)/2}\frac{x^k}{(1+x)_q^{n+k}},
\end{equation}
and
\begin{equation}
s^q_{n,k}(t)= E(-[n]_qt)\frac{([n]_qt)^k}{[k]_q!}.
\end{equation}
In case $q = 1,$ the above operators reduce to the Baskakov--Sz\'{a}sz operators \cite{gs}. \\
\indent Later, Mishra and Sharma \cite{MS} introduced a new Stancu type generalization of $q$-Baskakov-Sz\'{a}sz operators is defined as


\begin{equation}
\label{eq2}
\mathfrak{D}_{n}^{(\alpha,\beta)}(f;q;x)= [n]_q \sum_{k=0}^{\infty}p^q_{n,k}(x)\int_{0}^{q/1-q^n}q^{-k-1}s^q_{n,k}(t)f\bigg(\frac{[n]_qtq^{-k}+\alpha}{[n]_q+\beta}\bigg) d_qt ,
\end{equation}
where $p^q_{n,k}(x)$ and $s^q_{n,k}(t)$ are Baskakov and Sz\'{a}sz basis function respectively, defined as above. The operators $\mathcal{D}_{n}^{(\alpha,\beta)}(f,q;x)$ in (\ref{eq2}) are called $q$-Baskakov-Sz\'{a}sz-Stancu operators. For $\alpha=0,\,\beta=0$ the operators (\ref{eq2}) reduce to the operators (\ref{eq1}).\\

In the recent years, Stancu generalization of the certain operators introduced by several researchers and obtained different type of approximation properties of many operators, we refer some of the important
papers in this direction as \cite{DK,RN,
KJ2}. Recently, Mishra et al. (\cite{V14}, \cite{V2014}) have established very
interesting results on approximation properties of various functional classes
using different types of positive linear summability operators.

Before proceeding further, let us give some basic definitions and notations from
$q$-calculus. Such notations can be found in (\cite{er}, \cite{kc}). We consider $q$ as a real number satisfying $0<q<1$.\\
For

\begin{equation*} \displaystyle [n]_q = \left\{ \begin{array}{ll} \frac{1-q^n}{1-q}, & \hbox{$q\neq1$}, \\
n,& \hbox{$q=1$}, \end{array} \right.
 \end{equation*}
and

\begin{eqnarray*}
[n]_q!=\left\{
\begin{array}{ll}
[n]_q[n-1]_q[n-2]_q...[1]_q,  & \hbox{$n=1,2,...$},\\
1,& \hbox{$n=0$.}
\end{array}\right.
\end{eqnarray*}
Then for $q >0$ and integers $n, k, k \geq n \geq 0$, we have\\
$$[n+1]_{q}=1+q[n]_q  ~~~~{ \text {  and }}~~~ [n]_q+q^n[k-n]_q=[k]_q.$$
We observe that

\begin{equation*} (1+x)_q^n=(-x;q)_n= \left\{ \begin{array}{ll} (1+x)(1+qx)(1+q^2x)...(1+q^{n-1}x), & \hbox{$n=1,2,...$},\\
1,& \hbox{$n=0.$} \end{array}\right.\end{equation*}
Also, for any real number $\alpha$, we have

$$(1+x)_q^\alpha =\frac{(1+x)_q^\infty}{(1+q^\alpha x)_q^\infty}.$$
In special case, when $\alpha$ is a whole number, this definition coincides with the above definition.\\
The $q$-Jackson integral and $q$-improper integral defined as

\begin{equation*}\int_0^a f(x) d_qx=(1-q)a\sum_{n=0}^\infty f(aq^n)q^n\end{equation*}
and
\begin{equation*}\int_0^{\infty/A} f(x) d_qx=(1-q)a\sum_{n=0}^\infty f\left(\frac{q^n}{A}\right)\frac{q^n}{A},\end{equation*}
provided sum converges absolutely.\\
The $q$-analogues of the exponential function $e^x$ (see \cite{kc}), used here is defined as
\begin{equation*} E_q(z)=\prod_{j=0}^{\infty}(1+(1-q)q^jz)=\sum_{k=0}^{\infty}q^{k(k-1)/2} \frac{z^k}{[k]_q!}=(1+(1-q)z)_q^\infty, \ |q|<1,
\end{equation*}
where $(1-x)_q^{\infty}=\prod_{j=0}^{\infty}(1-q^jx).$\\

\section{Moment estimates}
\begin{lem}\label{l1}\cite{vj}: The following hold:
\begin{enumerate}
\item $\mathcal{D}_n(1,q;x)=1,$
\item $\mathcal{D}_n(t,q;x)=x+\frac{q}{[n]_q},$
\item $\displaystyle \mathcal{D}_n(t^2,q;x)={\bigg(1+\frac{1}{q[n]_q}\bigg)}x^2+\frac{x}{[n]_q}(1+q(q+2))+
    \frac{q^2(1+q)}{[n]_q^2}.$
\end{enumerate}
\end{lem}

\begin{lem}\label{l2}\cite{MS} The following  hold:
\begin{enumerate}
\item $\mathfrak{D}_{n}^{(\alpha,\beta)}(1;q;x)=1,$
\item $\displaystyle\mathfrak{D}_{n}^{(\alpha,\beta)}(t;q;x)=\frac{[n]_qx+q+\alpha}{[n]_q+\beta},$
\item $\displaystyle \mathfrak{D}_{n}^{(\alpha,\beta)}(t^2;q;x) =\bigg(\frac{[n]_q(q[n]_q+1)}{q([n]_q + \beta)^2}\bigg)x^2
    + \left(\frac{(1+q(q+2))[n]_q+2\alpha [n]_q}{([n]_q+\beta)^2}\right)x + \frac{q^2(1+q)+2q\alpha+\alpha^2}{([n]_q +\beta)^2}$.
\end{enumerate}
\end{lem}

\section{Korovkin type statistical approximation properties}
The idea of statistical convergence was introduced independently by Steinhaus
\cite{ST}, Fast \cite{HF} and Schoenberg \cite{98}.
In approximation theory, the concept of statistical convergence was used in the year 2002 by Gadjiev and
Orhan \cite{GO}. They proved the Bohman-Korovkin type approximation theorem
for statistical convergence. It was shown that the statistical versions are stronger than the classical ones.

Korovkin type approximation theory has also many useful connections, other than classical approximation theory, in other branches of mathematics (see Altomare and Campiti in \cite{ac22}).\\
 ~\indent Now, we recall the concept of statistical convergence for
sequences of real numbers which was introduced by Fast \cite{HF} and Mishra
et al. \cite{KJ}.

\parindent=8mm Let $K\subseteq \mathbb{N}$ and $K_{n}=\left\{ j\leq n:j\in
K\right\}.$ Then the $natural~density$ of $K$ is defined by $\delta
(K)={\lim\limits_{n}}~ n^{-1}|K_{n}|$ if the limit exists, where $|K_{n}|$ denotes the
cardinality of the set $K_{n}$.

\parindent=8mm A sequence $x=(x_{j})_{j\geq1}$ of real numbers is said to be $%
statistically$ $convergent$ to $L$ provided that for every $\epsilon >0$ the
set $\{j\in \mathbb{N}:|x_{j}-L|\geq \epsilon \}$ has natural density zero,
i.e. for each $\epsilon >0$,
\begin{equation*}
\lim\limits_{n}\frac{1}{n}|\{j\leq n:|x_{j}-L|\geq \epsilon \}|=0.
\end{equation*}
It is denoted by $st-\lim\limits_{n}x_{n}=L$.\\

In \cite{DO} DoÄŸru and Kanat defined the Kantorovich-type modification of Lupa\c{s} operators as follows:

\begin{equation}\label{dk}
\tilde{R}_n(f;q;x)=[n+1]\sum_{k=0}^{n}\bigg(\int_{\frac{[k]}{[n+1]}}^{\frac{[k+1]}{[n+1]}} ~ f(t) d_{q}t\bigg)\left(\begin{array}{c}n \\k \end{array}\right) \frac{q^{-k}q^{k(k-1)/2} x^k (1-x)^{(n-k)}}{(1-x+qx)\cdots(1-x+q^{n-1} x)}.
\end{equation}
DoÄŸru and Kanat \cite{DO} proved the following statistical Korovkin-type approximation
theorem for operators (\ref{dk}).

\begin{theorem}
Let $q:=(q_n),~ 0 < q < 1$, be a sequence satisfying the following conditions:\\
\begin{equation}\label{thm1}
 st-\lim_n q_n = 1,~ st - \lim_n q_n^n = a ~( a < 1)~ and~ st - \lim_n \frac{1}{[n]_q} = 0,
\end{equation}
then if $f$ is any monotone increasing function defined on $[0, 1]$, for the positive linear operators $\tilde{R}_n(f;q;x)$, then
$$ st- \lim_n {\Vert \tilde{R}_n(f;q;\cdot) - f \Vert}_{C[0, 1]} = 0$$ holds.
\end{theorem}
In \cite{do26} Do\u{g}ru gave some examples so that $(q_{n})$ is statistically convergent to $1$ but it may not convergent to $1$ in the ordinary case.\\
\newline
Now, we consider a sequence $q=(q_{n}),$ $q_{n}\in $ $(0,1),$ such that

\begin{equation}\label{a1}
\lim\limits_{n\rightarrow \infty }q_{n}=1.
\end{equation}
The condition (\ref{a1}) guarantees that
$[n]_{q_{n}}\rightarrow
\infty $ as $n\rightarrow \infty .$

\begin{theorem} \label{thm2}
Let ${{\mathfrak{D}^{(\alpha , \beta)}_{n}}}$ be the sequence of the operators (\ref{eq2}) and the sequence
$q = (q_n)$ satisfies (\ref{thm1}). Then for any function $f \in C[0,\nu] \subset C[0,\infty),~ \nu > 0 ,$ we have
\begin{equation}
st- \lim_n \Vert {{\mathfrak{D}^{(\alpha , \beta)}_n}}( f; q; \cdot) - f \vert = 0,
\end{equation}
where $C[0,\nu]$ denotes the space of all real bounded functions $f$ which are continuous in
$[0,\nu].$
\end{theorem}
\hspace{-0.9cm}
\textbf{Proof.}
Let $f_{i}= t^i,$ where $i=0,1,2.$   Using ${\mathfrak{D}^{(\alpha, \beta)}_n}(1;q_n;x)=1,$ it is clear that\\
$$st-\lim_{n}\|{\mathfrak{D}^{(\alpha, \beta)}_n}(1;q_n;x)-1\|=0.$$
Now by Lemma (\ref{l2})(ii), we have

\begin{equation*}
\lim_{n\rightarrow \infty }\Vert \mathfrak{D}_{n}^{(\alpha,\beta)}(t;q_{n};x)-x\Vert
=\left\|\frac{[n]_qx+q+\alpha}{[n]_q+\beta}-x\right\|\leq \frac{(q+\alpha)}{[n]_q + \beta} x + \frac{\beta}{([n]_q + \beta)}.
\end{equation*}
For given $\epsilon >0$, we define the following sets:
\begin{equation*}
 L=\{k:\Vert{\mathfrak{D}^{(\alpha, \beta)}_n}(t;q_{k};x)-x\Vert\geq \epsilon\},
\end{equation*}
and

\begin{equation}\label{e}
L'=\left\{k: \frac{(q+\alpha)}{[k]_q+\beta}x+\frac{\beta}{[k]_q + \beta}\geq \epsilon\right\}.
\end{equation}
It is obvious that $L\subset L^{\prime }$, it can be written as

\begin{equation*}
\delta \left( \{k\leq n:\Vert \mathfrak{D}_{n}^{(\alpha,\beta)}(t;q_{k};x)-x\Vert \geq
\epsilon \}\right) \leq \delta \left( \{k \leq n: ~ \frac{(q+\alpha)}{[k]_q+\beta}x + \frac{\beta}{[k]_q + \beta}\geq \epsilon \}\right).
\end{equation*}
By using (\ref{thm1}), we get

$$st-\lim_{n}\left(\frac{(q+\alpha)}{[n]_q+\beta}x + \frac{\beta}{[n]_q + \beta}\right)=0.$$
So, we have

$$\delta \bigg(\big\{k \leq n: ~ \frac{(q+\alpha)}{[n]_q+\beta}x + \frac{\beta}{[n]_q + \beta}\geq \epsilon \big\}\bigg) = 0,$$
then
$$st-\lim_{n}\|{\mathfrak{D}^{(\alpha, \beta)}_n}(t;q_n;x)-x\|=0.$$
Similarly, by Lemma (\ref{l2})(iii), we have

\begin{align*}
&\Vert \mathfrak{D}_{n}^{(\alpha,\beta)}(t^2;q_{n};x)-x^2\Vert\\& = \left\|\bigg(\frac{[n]_q(q[n]_q+1)}{q([n]_q + \beta)^2}\bigg)x^2
    + \left(\frac{(1+q(q+2))[n]_q+2\alpha [n]_q}{([n]_q+\beta)^2}\right)x+ \frac{q^2(1+q)+2q\alpha+\alpha^2}{([n]_q +\beta)^2}
- x^2\right\|\\
&\leq \bigg|{ \bigg(\frac{[n]_q(q[n]_q+1)}{q([n]_q + \beta)^2}\bigg)-1}\bigg| \nu^2+ \bigg| \left(\frac{(1+q(q+2))[n]_q+2\alpha [n]_q}{([n]_q+\beta)^2}\right)\bigg| \nu + \bigg|\frac{q^2(1+q)+2q\alpha+\alpha^2}{([n]_q +\beta)^2}\bigg|\\
&\leq {\mu}^2 \left[{ \bigg(\frac{[n]_q(q[n]_q+1)}{q([n]_q + \beta)^2}-1 \bigg)} + \left(\frac{(1+q(q+2))[n]_q+2\alpha [n]_q}{([n]_q+\beta)^2}\right) + \bigg(\frac{q^2(1+q)+2q\alpha+\alpha^2}{([n]_q +\beta)^2}\bigg)\right]
\end{align*}
where $\mu^2 = \max\{\nu^2,\nu,1\}=\nu^2.$\\
Now, if we choose

\begin{equation*}
\alpha _{n}=  \bigg(\frac{[n]_q(q[n]_q+1)}{q([n]_q + \beta)^2}-1 \bigg),
\end{equation*}%
\begin{equation*}
\beta _{n}=\left(\frac{(1+q(q+2))[n]_q+2\alpha [n]_q}{([n]_q+\beta)^2}\right),
\end{equation*}%
\begin{equation*}
\gamma _{n}= \bigg(\frac{q^2(1+q)+2q\alpha+\alpha^2}{([n]_q +\beta)^2}\bigg),
\end{equation*}%
now using $(\ref{thm1}),$ we can write
\begin{equation} \label{2.1}
st-\lim_{n\rightarrow \infty }\alpha _{n}=0=st-\lim_{n\rightarrow \infty
}\beta _{n}=st-\lim_{n\rightarrow \infty }\gamma _{n}.
\end{equation}
Now for given $\epsilon >0$, we define the following four sets
\begin{equation*}
\mathcal{U}=\{k: \Vert {\mathfrak{D}^{(\alpha, \beta)}_n}(t^{2};q_{k};x)-x^{2}\Vert \geq \epsilon \},
\end{equation*}%
\begin{equation*}
\mathcal{U}_{1}=\{k:\alpha _{k}\geq \frac{\epsilon }{{\mu}^2}\},
\end{equation*}%
\begin{equation*}
\mathcal{U}_{2}=\{k:\beta _{k}\geq \frac{\epsilon }{{\mu}^2}\},
\end{equation*}%

\begin{equation*}
\mathcal{U}_{3}=\{k:\gamma _{k}\geq \frac{\epsilon }{{\mu}^2}\}.
\end{equation*}%
It is obvious that $\mathcal{U}\subseteq {\mathcal{U}_{1}\cup \mathcal{U}_{2}\cup \mathcal{U}_{3}}$. Then, we obtain%

\begin{eqnarray*}
\delta \big(\{k &\leq & n:\Vert \mathfrak{D}_{n}^{(\alpha,\beta)}(t^{2};q_{n};x)-x^{2}\Vert \geq
\epsilon \}\big) \\&&~~~~~ \leq \delta \big(\{k\leq n:\alpha _{k}\geq \frac{\epsilon }{{\mu}^2}\}\big) + \delta
\big(\{k\leq n:\beta _{k}\geq \frac{\epsilon }{{\mu}^2}\}\big)+\delta \big(\{k\leq n:\gamma
_{k}\geq \frac{\epsilon }{{\mu}^2}\}\big).
\end{eqnarray*}%
\newline
Using (\ref{2.1}), we get%

\begin{equation*}
st-\lim_{n \rightarrow \infty }\Vert  {\mathfrak{D}^{(\alpha, \beta)}_n}(t^{2};q_{n};x)-x^{2}\Vert =0.
\end{equation*}
Since,
\begin{eqnarray*}
\Vert  {\mathfrak{D}^{(\alpha, \beta)}_n}(f;q_{n};x)-f\Vert \leq \Vert {\mathfrak{D}^{(\alpha, \beta)}_n}(t^{2};q_{n};x)-x^{2}\Vert +\Vert  {\mathfrak{D}^{(\alpha, \beta)}_n}
(t;q_{n};x)-x\Vert +\Vert {\mathfrak{D}^{(\alpha, \beta)}_n}(1;q_{n};x)-1\Vert ,\newline
\end{eqnarray*}%
we get
\begin{eqnarray*}
st-\lim_{n\rightarrow \infty }\Vert {\mathfrak{D}^{(\alpha, \beta)}_n}(f;q_{n};x)-f\Vert
&\leq &  st-\lim_{n\rightarrow \infty }\Vert {\mathfrak{D}^{(\alpha, \beta)}_n}(t^{2};q_{n};x)-x^{2}\Vert\\
&&+~st-\lim_{n\rightarrow \infty }\Vert{\mathfrak{D}^{(\alpha, \beta)}_n}(t;q_{n};x)-x\Vert\\
&&+~st-\lim_{n\rightarrow \infty }\Vert {\mathfrak{D}^{(\alpha, \beta)}_n}(1;q_{n};x)-1\Vert ,
\end{eqnarray*}
which implies that
\begin{equation*}
st-\lim_{n\rightarrow \infty }\Vert  {\mathfrak{D}^{(\alpha, \beta)}_n}(f;q_{n};x)-f\Vert =0.
\end{equation*}
This completes the proof of theorem. \qed

\section{Weighted statistical approximation}
In this section, we obtain the Korovkin type weighted
statistical approximation by the operators defined in (\ref{eq2}).
A real function $\rho $ is called a weight function if it is continuous on $%
\mathbb{R}$ and $\lim\limits_{\mid x\mid \rightarrow \infty }\rho (x)=\infty
,~\rho (x)\geq 1$ for all $x\in \mathbb{R}$.

Let by $B_{\rho }(\mathbb{R})$  denote the weighted space of real-valued
functions $f$ defined on $\mathbb{R}$ with the property $\mid f(x)\mid \leq
M_{f}~\rho (x)$ for all $x\in \mathbb{R}$, where $M_{f}$ is a constant
depending on the function $f$. We also consider the weighted subspace $%
C_{\rho }(\mathbb{R})$ of $B_{\rho }(\mathbb{R})$ given by $C_{\rho }(%
\mathbb{R})=\{f\in B_{\rho }(\mathbb{R}){:}$ $f$ continuous on $\mathbb{R}%
\} $. Note that $B_{\rho }(\mathbb{R})$ and $C_{\rho }(\mathbb{R})$ are
Banach spaces with $\Vert f\Vert _{\rho }=\sup\limits_{x\in R}\frac{\mid
f(x)\mid }{\rho (x)}.$ In case of weight function $~\rho _{0}=1+x^{2},$ we have $%
\Vert f\Vert _{\rho _{0}}=\sup\limits_{x\in R}\dfrac{\mid f(x)\mid }{1+x^{2}}%
.$

Now we are ready to prove our main result as follows:\newline
\begin{theorem}\label{thm3}
Let $\mathfrak{D}_{n}^{(\alpha,\beta)}$ be the sequence of the operators (\ref{eq2})and the sequence $q=(q_{n})$ satisfies (\ref{thm1}). Then for any function $f\in C_{B}[0,\infty ),$ we have

\begin{equation*}
st-\lim_{n\rightarrow \infty }{\Vert }\mathfrak{D}_{n}^{(\alpha,\beta)}(f;q_{n};.)-f{\Vert }%
_{\rho _{0}}=0.
\end{equation*}
\end{theorem}
\hspace{-0.9cm}
\textbf{Proof.} By Lemma (\ref{l2})(iii), we have  ${\mathfrak{D}^{(\alpha, \beta)}_n}(t^2; q_n; x) \leq Cx^2,$ where $C$ is a positive constant, ${\mathfrak{D}^{(\alpha, \beta)}_n}(f; q_n ; x)$ is a sequence of positive linear operators acting from $C_{\rho}[0,\infty)$ to $B_{\rho}[0,\infty)$. \\
Using ${\mathfrak{D}^{(\alpha, \beta)}_n}(1; q_n; x) = 1,$ it is clear that\\

$$st-\lim_{n}\|{\mathfrak{D}^{(\alpha, \beta)}_n}(1;q_{n};x)-1\|_{\rho_{0}}=0.$$
Now, by Lemma (\ref{l2})(ii), we have

$$\|{\mathfrak{D}^{(\alpha, \beta)}_n}(t; q_n; x)-x\|_{\rho_{0}}= \sup_{x\in [0,\infty)} \frac{|{\mathfrak{D}^{(\alpha, \beta)}_n}(t;q_n;x)-x|}{1+x^2} \leq \frac{(q+\alpha)}{[n]_q + \beta}  + \frac{\beta}{([n]_q + \beta)}.$$
Using (\ref{thm1}), we get

 $$st-\lim_{n}\left( \frac{(q+\alpha)}{[n]_q + \beta} + \frac{\beta}{([n]_q + \beta)}\right)=0,$$
then $$st-\lim_{n}\|{\mathfrak{D}^{(\alpha, \beta)}_n}(t;q_n;x)-x\|_{\rho_{0}}=0.$$
Finally, by Lemma(\ref{l2})(iii), we have

\begin{align*}
\Vert {\mathfrak{D}^{(\alpha, \beta)}_n}(t^2;q_{n};x)-x^2\Vert_{\rho_{0}}
&\leq \left({\frac{[n]_q(q[n]_q+1)}{q([n]_q + \beta)^2}-1} \right)\sup_{x\in[0,\infty)}\frac{x^2}{1+x^2}\\
&~~~+\left(\frac{(1+q(q+2))[n]_q+2\alpha [n]_q}{([n]_q+\beta)^2}\right)\sup_{x\in[0,\infty)}\frac{x}{1+x^2} + \frac{q^2(1+q)+2q\alpha+\alpha^2}{([n]_q +\beta)^2}\\
& \leq \left({\frac{[n]_q(q[n]_q+1)}{q([n]_q + \beta)^2}-1} \right)+\left(\frac{(1+q(q+2))[n]_q+2\alpha [n]_q}{([n]_q+\beta)^2}\right)+ \frac{q^2(1+q)+2q\alpha+\alpha^2}{([n]_q +\beta)^2}.\\
\end{align*}
Now, If we choose

\begin{equation*}
\alpha _{n}=  \bigg(\frac{[n]_q(q[n]_q+1)}{q([n]_q + \beta)^2}-1 \bigg),
\end{equation*}%
\begin{equation*}
\beta _{n}=\left(\frac{(1+q(q+2))[n]_q+2\alpha [n]_q}{([n]_q+\beta)^2}\right),
\end{equation*}%
\begin{equation*}
\gamma _{n}= \bigg(\frac{q^2(1+q)+2q\alpha+\alpha^2}{([n]_q +\beta)^2}\bigg),
\end{equation*}%
then by $(\ref{thm1}),$ we can write

\begin{equation}\label{3.1}
st-\lim_{n\rightarrow \infty }\alpha _{n}=0=st-\lim_{n\rightarrow \infty
}\beta _{n}=st-\lim_{n\rightarrow \infty }\gamma _{n}.
\end{equation}
Now for given $\epsilon >0$, we define the following four sets:

\begin{equation*}
S=\{k: \Vert {\mathfrak{D}^{(\alpha, \beta)}_n}(t^{2};q_{k};x)-x^{2}\Vert_{\rho_{0}} \geq \epsilon \},
\end{equation*}%
\begin{equation*}
S_{1}=\lbrace k:\alpha _{k}\geq \frac{\epsilon}{3}\rbrace,
\end{equation*}%
\begin{equation*}
S_{2}=\{k:\beta _{k}\geq \frac{\epsilon }{3}\},
\end{equation*}%
\begin{equation*}
S_{3}=\{k:\gamma _{k}\geq \frac{\epsilon }{3}\}.
\end{equation*}%
It is obvious that $S  \subseteq {S_{1}\cup S_{2}\cup S_{3}}$. Then, we obtain%
\begin{eqnarray*}
\delta \big(\{k &\leq &n:\Vert {\mathfrak{D}^{(\alpha, \beta)}_n}(t^{2};q_{n};x)-x^{2}\Vert_{\rho_{0}} \geq
\epsilon \}\big) \\
&&~~~~~~\leq \delta\big(\{k\leq n:\alpha _{k}\geq \frac{\epsilon }{3}\}\big)+\delta
\big(\{k \leq n:\beta _{k}\geq \frac{\epsilon }{3}\}\big)+\delta\big(\{k\leq n:\gamma
_{k}\geq \frac{\epsilon }{3}\}\big).
\end{eqnarray*}%
\newline
Using (\ref{3.1}), we get%
\begin{equation*}
st-\lim_{n\rightarrow \infty }\Vert {\mathfrak{D}^{(\alpha, \beta)}_n}(t^{2};q_{n};x)-x^{2}\Vert_{\rho_{0}} =0.
\end{equation*}
Since
\begin{eqnarray*}
&&\Vert{\mathfrak{D}^{(\alpha, \beta)}_n}(f;q_{n};x)-f\Vert_{\rho_{0}}\\
&&~~~~~~~~~\leq \Vert {\mathfrak{D}^{(\alpha, \beta)}_n}(t^{2};q_{n};x)-x^{2}\Vert_{\rho_{0}} +\Vert  {\mathfrak{D}^{(\alpha, \beta)}_n}
(t;q_{n};x)-x\Vert_{\rho_{0}}+\Vert {\mathfrak{D}^{(\alpha, \beta)}_n}(1;q_{n};x)-1\Vert_{\rho_{0}} ,\newline
\end{eqnarray*}%
we get
\begin{eqnarray*}
st-\lim_{n\rightarrow \infty }\Vert{\mathfrak{D}^{(\alpha, \beta)}_n}(f;q_{n};x)-f\Vert_{\rho_{0}}
&\leq &st-\lim_{n\rightarrow \infty }\Vert{\mathfrak{D}^{(\alpha, \beta)}_n}(t^{2};q_{n};x)-x^{2}\Vert_{\rho_{0}}\\&& +~st-\lim_{n\rightarrow \infty }\Vert{\mathfrak{D}^{(\alpha, \beta)}_n}(t;q_{n};x)-x\Vert_{\rho_{0}}\\
&& + ~st-\lim_{n\rightarrow \infty }\Vert  {\mathfrak{D}^{(\alpha, \beta)}_n}(1;q_{n};x)-1\Vert_{\rho_{0}} ,
\end{eqnarray*}%
which implies that
\begin{equation*}
st-\lim_{n\rightarrow \infty }\Vert{\mathfrak{D}^{(\alpha, \beta)}_n}(f;q_{n};x)-f\Vert_{\rho_{0}} =0.
\end{equation*}
This completes the proof of the theorem.\qed


\section{Rates of statistical convergence}
In this section, by using the modulus of continuity, we will study rates of statistical convergence of operators (\ref{eq2}) and Lipschitz type maximal functions are introduced.

\begin{lem}\label{l3}
Let $0<q<1$ and $a\in[0,bq],~ b>0.$ The inequality
\begin{equation}
\int_{a}^{b}|t-x|d_qt\leq \left( \int_{a}^{b}|t-x|^2d_qt\right)^{1/2} \left(\int_{a}^{b} d_qt \right)^{1/2}
\end{equation}
is satisfied.
\end{lem}
Let $C_B[0,\infty),$ the space of all bounded and continuous functions on $[0,\infty)$ and $x \geq 0.$
Then, for $\delta>0,$ the modulus of continuity of $f$ denoted by $\omega(f;\delta )$ is defined to be
$$\omega(f;\delta)= \sup_{|{t- x}|\leq {\delta}}|f(t)-f(x)|,~t\in[0,\infty).$$
It is known that $\lim\limits_{\delta\rightarrow 0}\omega(f ; \delta) = 0$ for $f\in C_B[0,\infty)$ and also, for any $\delta > 0$ and each $t,~x\geq 0,$ we have

\begin{equation}\label{2.2}
|f(t)-f(x)|\leq \omega(f;\delta)\left(1+\frac{|t-x|}{\delta}\right).
\end{equation}

\begin{theorem}
Let $(q_n)$ be a sequence satisfying (\ref{thm1}). For every non-decreasing $f\in C_B[0,\infty), ~x\geq 0$ and $n\in \mathbb{N},$ we have
$$|{\mathfrak{D}^{(\alpha, \beta)}_n}(f;q_{n};x)-f(x)| \leq 2\omega(f;\sqrt{\delta_n(x)}),$$
where

\begin{eqnarray*}\label{d1}
\delta_n{(x)}&=&\left(\frac{[n]_q(q[n]_q+1)}{q([n]_q+\beta)^2}+1-\frac{2[n]_q}{[n]_q+\beta}\right)x^2\\
&&+\left(\frac{[n]_q+q^2[n]_q-2\alpha\beta-2q\beta}{([n]_q+\beta)^2}\right)x+\frac{q^2(1+q)+2q\alpha+\alpha^2}{([n]_q+\beta)^2}.
\end{eqnarray*}
\end{theorem}

\hspace{-0.9cm}
\textbf{Proof.}
Let  $f\in C_B[0,\infty)$ be a non-decreasing function and $x\geq 0$. Using linearity and positivity of the operators
 ${\mathfrak{D}}^{(\alpha, \beta)}_n$ and then applying (\ref{2.2}), we get for $\delta > 0$
\begin{eqnarray*}
|{\mathfrak{D}^{(\alpha, \beta)}_n}(f;q_{n};x)-f(x)| &\leq &{\mathfrak{D}^{(\alpha, \beta)}_n}\big(|f(t)-f(x)|;q_{n};x\big)\\
&\leq & \omega(f,\delta)\big\{{\mathfrak{D}^{(\alpha, \beta)}_n}(1;q_{n};x)+
\frac{1}{\delta}{\mathfrak{D}^{(\alpha, \beta)}_n}(|t-x|;q_{n};x) \big\}.
\end{eqnarray*}
Taking
${\mathfrak{D}^{(\alpha, \beta)}_n}(1;q_{n};x) =1$ and using Cauchy-Schwartz inequality,  we have

\begin{align*}
|{\mathfrak{D}^{(\alpha, \beta)}_n}(f;q_{n};x)-f(x)|
&\leq  \omega(f;\delta)\bigg\{1+\frac{1}{\delta}{\bigg({\mathfrak{D}^{(\alpha, \beta)}_n}((t-x)^2;q_{n};x\big)}^{1/2}
{{\mathfrak{D}^{(\alpha, \beta)}_n}(1;q_{n};x)}^{1/2}\bigg)\bigg\}
\end{align*}

\begin{align*}
&\leq  \omega(f;\delta)\bigg[1+\frac{1}{\delta}\bigg\{\left(\frac{[n]_q(q[n]_q+1)}{q([n]_q+\beta)^2}+1-\frac{2[n]_q}{[n]_q+\beta}\right)x^2\\
&~~~~~~~~~~~+\left(\frac{[n]_q+q^2[n]_q-2\alpha\beta-2q\beta}{([n]_q+\beta)^2}\right)x+\frac{q^2(1+q)+2q\alpha+\alpha^2}{([n]_q+\beta)^2}\bigg\}^{1/2}\bigg].
\end{align*}
Taking $q = (q_n),$ a sequence satisfying (\ref{thm1}) and choosing $\delta = \delta_n(x)$ as in (\ref{d1}), the theorem is proved.\qed \\


Now we will give an estimate concerning the rate of approximation by means
of Lipschitz type maximal functions.\\
 In \cite{L}, Lenze introduced a Lipschitz type maximal function as%
\begin{equation*}
{f}_{\alpha }(x,y)=\sup\limits_{t>0,t\neq x}\frac{\mid f(t)-f(x)\mid }{{%
\mid }t-x{\mid }^{\alpha }}.
\end{equation*}%
\newline
In \cite{adb24}, the Lipschitz type maximal function space on $E\subset \lbrack
0,\infty )$ is defined as follows%
\begin{equation*}
\tilde{V}_{\alpha,E }=\{f=\sup (1+x)^{\alpha }~ {f}_{\alpha }(x,y)\leq M%
\frac{1}{(1+y)^{\alpha }};x\geq 0~{and}~y~\in E\},
\end{equation*}%
where $f$ is bounded and continuous function on $[0,\infty )$, $M$ is a
positive constant and $0<\alpha \leq 1$. \newline
Also, let $d(x,E)$ be the distance between $x$ and $E,$ that is,
$$d(x,E)= \inf \{|x-y|; y \in E\}.$$

\begin{theorem}\label{thm4}
 If ${\mathfrak{D}^{(\alpha, \beta)}_n}$ be defined by (\ref{eq2}), then for all $%
f\in \tilde{V}_{\alpha ,E}$
\begin{equation}\label{e2.4.3}
\mid {\mathfrak{D}^{(\alpha, \beta)}_n}(f;q_{n};x)-f(x)\mid \leq M(\delta _{n}^{\frac{\alpha }{2}%
} ~ + ~d(x,E)),
\end{equation}%
where
\begin{equation}\label{d2}
\delta_n{(x)}=\left(\frac{[n]_q(q[n]_q+1)}{q([n]_q+\beta)^2}+1-\frac{2[n]_q}{[n]_q+\beta}\right)x^2
+\left(\frac{[n]_q+q^2[n]_q-2\alpha\beta-2q\beta}{([n]_q+\beta)^2}\right)x+\frac{q^2(1+q)+2q\alpha+\alpha^2}{([n]_q+\beta)^2}.
\end{equation}%
\end{theorem}
\hspace{-0.9cm}
\textbf{Proof.} Let $x_0\in \bar{E}$, where $\bar{E}$ denote the closure of the set $E$. Then
we have\newline
\begin{equation*}
\mid f(t)-f(x)\mid \leq \mid f(t)-f(x_{0})\mid +\mid f(x_{0})-f(x)\mid .
\end{equation*}%
Since ${\mathfrak{D}^{(\alpha, \beta)}_n}$ is a positive and linear operators, $f\in \tilde{V}%
_{\alpha ,E}$ and using the above inequality

\begin{equation*}
\mid {\mathfrak{D}^{(\alpha, \beta)}_n}(f;q_{n};x)-f(x)\mid ~ \leq {\mathfrak{D}^{(\alpha, \beta)}_n}(\mid
f(t)-f(x_{0})\mid ;q_{n};x) + (\mid f(x_{0})-f(x)\mid){\mathfrak{D}^{(\alpha, \beta)}_n}(1;q_{n};x)
\end{equation*}
\begin{equation}\label{e2.4.5}
\leq M\left( {\mathfrak{D}^{(\alpha, \beta)}_n}({\mid }t-x_{0}{\mid }^{\alpha };q_{n};x)+{\mid }%
x-x_{0}{\mid }^{\alpha }{\mathfrak{D}^{(\alpha, \beta)}_n}(1;q_{n};x)\right) .
\end{equation}%
Therefore, we have
\begin{equation*}
{\mathfrak{D}^{(\alpha, \beta)}_n}\left( {\mid }t-x_{0}{\mid }^{\alpha };q_{n};x\right) \leq
{\mathfrak{D}^{(\alpha, \beta)}_n}({\mid }t-x{\mid }^{\alpha };q_{n};x)+{\mid }x-x_{0}{\mid }%
^{\alpha }{\mathfrak{D}^{(\alpha, \beta)}_n}(1;q_{n};x).
\end{equation*}%
Now, we take $p=\frac{2}{\alpha }$ and $q=\frac{2}{(2-\alpha) }$ and by using the H\"{o}lder's inequality, one can write

\begin{equation*}
{\mathfrak{D}^{(\alpha, \beta)}_n}\left( {(}t-x)^{\alpha };q_{n};x\right) \leq {\mathfrak{D}^{(\alpha, \beta)}_n}\left( {(}t-x{)}^{2};q_{n};x\right) ^{\alpha/2}({\mathfrak{D}^{(\alpha, \beta)}_n}(1;q_{n};x){)}%
^{{(2-\alpha)}/{2}}
\end{equation*}%
\begin{equation*}
+{\mid }x-x_{0}{\mid }^{\alpha }{\mathfrak{D}^{(\alpha, \beta)}_n}(1;q_{n};x)
\end{equation*}
\begin{equation*}
=\delta _{n}^{\frac{\alpha }{2}} + {\mid }x-x_{0}{%
\mid }^{\alpha }.
\end{equation*}%
Substituting this in ($\ref{e2.4.5}$), we get ($\ref{e2.4.3}$).\newline
This completes the proof of the theorem.\newline \qed

\begin{remark}
 Observe that by the conditions in (\ref{thm1}), $$st-\lim\limits_n \delta_n = 0.$$\\ By (\ref{2.2}), we have
$$st-\lim\limits_n\omega(f;\delta_n) = 0.$$
This gives us the pointwise rate of statistical convergence of the operators  ${\mathfrak{D}^{(\alpha, \beta)}_n}(f;q_{n};x)$ to $f(x).$
\end{remark}

\begin{remark}
If we take $E=[0,\infty )$ in Theorem \ref{thm4} , since $d(x,E)=0,$
then we get for every $f\in \tilde{V}_{\alpha ,[0,\infty )}$
\begin{equation*}
{\mid }{\mathfrak{D}^{(\alpha, \beta)}_n}(f;q_{n};x)-f(x){\mid }\leq M\delta _{n}^{\frac{%
\alpha }{2}}
\end{equation*}%
where $\delta _{n}$ is defined as in (\ref{d2}).
\end{remark}

\begin{remark}
 By using (\ref{3.1}), It is easy to verify that
\begin{equation*}
st-\lim_{n\rightarrow \infty }\delta _{n}=0.
\end{equation*}%
That is, the rate of statistical convergence of operators (\ref{eq2}) to the
function $f$ are estimated by means of Lipschitz type maximal functions.
\end{remark}

\vspace{0.2cm}
\hspace{-0.9cm}
\textsc{\bf Acknowledgment}\\
The authors would like to express their deep gratitude to the anonymous learned referee(s) and the editor for their valuable suggestions and constructive comments, which resulted in the  subsequent improvement of this research article. Special thanks are due to our great Masters and friend academicians Prof. Abdel-shafi Obada, Editor in Chief, Prof. A.M. El-Sayed, Associate Editor of JOEMS for kind cooperation, smooth behavior during communication and for their efforts to send the reports of the manuscript timely. The second and third authors P. Sharma and L.N. Mishra acknowledge the MHRD, New Delhi, India for supporting this research article. The authors declare that there is no conflict of interests regarding the publication of this research article.

\newpage
\hspace{-0.9cm}

\end{document}